\documentclass[12pt]{article}
\usepackage{amssymb}
\usepackage{graphicx}
\usepackage{amsfonts}
\usepackage{latexsym}
\usepackage{amsthm}
\usepackage{amsmath}

\usepackage{amssymb, amscd}

\usepackage{url}

\newtheorem{problem}{Problem}
\newtheorem{theo}[problem]{Theorem}

\newtheorem{defin}[problem]{Definition}
\newtheorem{prop}[problem]{Proposition}
\newtheorem{cor}[problem]{Corollary}
\newtheorem{lema}[problem]{Lemma}
\newtheorem{exam}[problem]{Example}

\begin{document}
\date{November 2007}
 \title{{ Cycle-free chessboard complexes\\ and symmetric homology of algebras}}

\author{Sini\v sa T.\ Vre\' cica\\ {\small Faculty of Mathematics}\\[-2mm] {\small Belgrade
University}\\[-2mm] {\small vrecica$@$matf.bg.ac.yu}
 \and Rade  T.\ \v Zivaljevi\' c\\ {\small Mathematical Institute}\\[-2mm] {\small SANU, Belgrade}\\[-2mm]
 {\small rade$@$mi.sanu.ac.yu} }

\maketitle
\begin{abstract}
Chessboard complexes and their relatives have been an important
recurring theme of topological combinatorics (see \cite{Ata04},
\cite{BLVZ}, \cite{FrHa98}, \cite{Ga79}, \cite{Jo-book},
\cite{Jo07}, \cite{ShaWa04}, \cite{Wa03}, \cite{Zie94},
\cite{ZV92}). Closely related ``cycle-free chessboard complexes''
have been recently introduced by Ault and Fiedorowicz in
\cite{AuFie07} and \cite{Fie07} as a tool for computing symmetric
analogues of the cyclic homology of algebras. We study
connectivity properties of these complexes and prove a result that
confirms a strengthened conjecture from \cite{AuFie07}.
\end{abstract}

\renewcommand{\thefootnote}{}
\footnotetext{The authors are supported by Grant 144026 of the
Ministry for Science of Serbia.}

\section{Introduction}
\label{sec:intro}

Chessboard complexes and their relatives are well studied objects
of topological combinatorics with applications in group theory,
representation theory, commutative algebra, Lie theory,
computational geometry, and combinatorics. The reader is referred
to \cite{Jo-book} and \cite{Wa03} for  surveys and to \cite{Jo07},
\cite{ShaWa04} for a guide to some of the latest developments.

Chessboard complexes originally appeared in \cite{Ga79} as coset
complexes of the symmetric group, closely related to Coxeter  and
Tits coset complexes. After that they have been rediscovered
several times. Among their avatars are ``complexes of partial
injective functions'' \cite{ZV92},  ``multiple deleted joins'' of
$0$-dimensional complexes  \cite{ZV92} (implicit in \cite{Sar91}),
the complex of all partial matchings in a complete bipartite
graph, the complex of all non-taking rook configurations
\cite{BLVZ} etc.

Recently a naturally defined subcomplex of the chessboard complex,
here referred to as ``cycle-free chessboard complex'', has emerged
in the context of stable homotopy theory (\cite{AuFie07} and
\cite{Fie07}). Ault and Fiedorowicz introduced this complex and
its suspension $Sym^{(p)}_*$ as a tool for evaluating the
symmetric analogue for the cyclic homology of algebras,
\cite{Lod98}. They conjectured that
$H_i\left(Sym^{(p)}_*\right)=0$ if $i<p/2$, and verified this
conjecture for some small values of $p$ and $i$.

In this paper we prove this conjecture (Theorem~\ref{thm:main}) by
showing that $Sym^{(p)}_*$ is actually $\gamma_p$-connected where
$\gamma_p = \left[\frac 23(p-1)\right]$ (see Corollary
\ref{posl}). We also show (Theorem~\ref{thm:tight1}) that this
result cannot be improved if $p=3k+2$ for some $k$ and give
evidence that the bound should be tight in the general case.

\subsection{Graph complexes}
\label{sec:graph complexes}

Chessboard complexes $\Delta_{m,n}$ and their relatives are
examples of {\em graph complexes}. A graph (digraph, multigraph)
complex is (in {\em topological combinatorics}) a family of graphs
(digraphs, multigraphs) on a given vertex set, closed under
deletion of edges. Monograph \cite{Jo-book}, based on the author's
Ph.D.\ degree thesis \cite{Jo-thesis}, serves as an excellent
source of information about applications of graph complexes in
algebraic and geometric/topological combinatorics and related
fields.

The appearance of a monograph solely devoted to the exposition and
classification of simplicial complexes of graphs is probably a
good sign of a relative maturity of the field. After decades of
development, some of the central research themes and associated
classes of examples have been isolated and explored, the technique
is codified and typical applications described.

However, the appearance of a relative of the chessboard complex in
the context of symmetric homology $HS_\ast(A)$ of algebras is
perhaps of somewhat non-standard nature and deserves a comment.

 Ault and Fiedorowicz showed in \cite{AuFie07} (Theorem 6) that
 there exists a spectral sequence converging strongly to
 $HS_\ast(A)$ with the $E^1$-term
\begin{equation}\label{eqn:spectral}
E^1_{p,q} = \bigoplus_{\overline{u}\in X^{p+1}/S_{p+1}}
\widetilde{H}_{p+q}(EG_{\overline{u}}\ltimes_{G_{\overline{u}}}
N\mathcal{S}_p/N\mathcal{S}_p'\, ;k).
\end{equation}
They emphasized (loc.\ cit.) the importance of the problem of
determining the homotopy type of the space
$N\mathcal{S}_p/N\mathcal{S}_p'$ and introduced a much more
economical complex $Sym_\ast^{(p)}$ which computes its homology.

The complex $Sym_\ast^{(p)}$ turned out to be isomorphic to the
suspension $\Sigma(\Omega_{p+1})$ of a subcomplex $\Omega_{p+1}$
of the chessboard complex $\Delta_{p+1}=\Delta_{p+1,p+1}$, one of
the well studied graph complexes!

\medskip
It is interesting to compare this development with the appearance
of {\em $2$-connected graph complexes} \cite{Vass-99} in the
computation of the $E^1$-term of the main Vassiliev spectral
sequence converging to the cohomology

\begin{equation}\label{eqn:knot}
H^i(\mathcal{K}\setminus\Sigma)\cong
\bar{H}_{\omega-i-1}(\Sigma)\cong \bar{H}_{\omega-i-1}(\sigma)
\end{equation}
of the space $\mathcal{K}\setminus\Sigma$ of non-singular knots in
$\mathbb{R}^n$. This spectral sequence arises from a filtration
$\sigma_1\supset\sigma_2\supset\ldots$ of a simplicial resolution
$\sigma$ of the space (discriminant) $\Sigma$ of singular knots in
$\mathcal{K}$. As a tool for computing
$E^1_{i,j}=\bar{H}_{i+j}(\sigma_i\setminus\sigma_{i-1})$,
Vassiliev \cite{Vass-99} introduced an auxiliary filtration of the
space $\sigma_i\setminus\sigma_{i-1}$. Complexes of $2$-connected
graphs naturally appear in the description of the $E^1$-term of
the spectral sequence associated to the auxiliary filtration.

It appears, at least on the formal level, that cycle-free
chessboard complexes $\Omega_n$ play the role, in the Ault and
Fiedorowicz approach to symmetric homology, analogous to the role
of $2$-connected graph complexes in Vassiliev's approach to the
homology of knot spaces.

\medskip
The homotopy type of the complex of (not) $2$-connected graphs was
(independently) determined by Babson, Bj\" orner, Linusson,
Shareshian, and Welker in \cite{BBLSW-99} and Turchin in
\cite{Tur-97}. This development stimulated further study of
connectivity graph properties (complexes), see chapter~VI of
\cite{Jo-book} (\cite{Jo-thesis}).

\section{Cycle-free chessboard complexes}

Chessboard complexes $\Delta_{m,n}$ are matching (graph) complexes
associated to complete bipartite graphs \cite{Jo-book},
\cite{ShaWa04}, \cite{Wa03}. However, they most naturally arise as
complexes of {\em admissible rook configurations} on  general
$m\times n$ chessboards.

A $(m\times n)$-chessboard is the set $A_{m,n}=[m]\times
[n]\subset \mathbb{Z}^2$ where (as usual in combinatorics)
$[n]=\{1,2,\ldots,n\}$. The associated {\em chessboard complex}
$\Delta(A_{m,n})=\Delta_{m,n}$ is defined as the (abstract,
simplicial) complex of all admissible or {\em non-taking rook
configurations} on the chessboard $A_{m,n}$. More generally, for
an arbitrary (finite) subset $A\subset \mathbb{Z}^2$, the
associated chessboard complex $\Delta(A)$ has $A$ for the set of
vertices and $S\in\Delta(A)$ if and only if for each pair $(i,j),
(i',j')$ of distinct vertices of $S$ both $i\neq i'$ and $j\neq
j'$. Also, we often denote by $\Delta_{X,Y}=\Delta(X\times Y)$ the
chessboard complex carried by the ``chessboard'' $X\times Y$ where
$X$ and $Y$ are not necessarily subsets of $\mathbb{Z}$.

Let $\Delta_n:=\Delta_{n,n}$ be the chessboard complex associated
to the canonical $(n\times n)$-chessboard $[n]\times [n]$,
similarly $\Delta_X:=\Delta_{X,X}$. Each top dimensional simplex
in $\Delta_n$ is essentially the graph $\Gamma_\phi:=\{(i,\phi(i))
\mid i\in [n]\}$ of a permutation $\phi : [n]\rightarrow [n]$. Any
other simplex $S\in\Delta_n$ arises as a top dimensional simplex
of the complex $\Delta(A\times B)$ where $A$ and $B$ are two
subsets of $[n]$ of equal size. Alternatively $S$ can be described
as the graph $\Gamma_\psi$ of a bijection $\psi : A\rightarrow B$,
which is sometimes referred to as a partial, injective function
(relation) defined on $[n]$.

\begin{defin}\label{def:cycle-free}
A non-taking rook configuration $S\subset [n]\times [n]$ of size
$n-1$ is {\em cycle-free} if there is a linear order  $\rho :
i_1\prec i_2\prec \ldots\prec i_n$ of elements of the set $[n]$
such that
$$
S = S_\rho = \{(i_1,i_2),(i_2,i_3),\ldots, (i_{n-1},i_{n})\}.
$$
Define $\Omega_n:= \bigcup_{\rho\in LO_n}~S_\rho\subset\Delta_n$
as the union, or $\subset$-ideal closure, of the collection of all
cycle-free configurations $S_\rho$, $\rho\in LO_n$, where $LO_n$
is the set of all linear orders on $[n]$.

Alternatively the complex $\Omega_n$ can be described as the
collection of all non-taking rook placements $S\in \Delta_n$ which
do not contain cycles, that is sub-configurations of the form
$\{(x_1,x_2),(x_2,x_3),\ldots,(x_m,x_1)\}$ for some $1\leq m\leq
n$. For this reason we call $\Omega_n$ the chessboard complex
without cycles or simply the {\em cycle-free} chessboard complex.
\end{defin}

In order to study functorial properties of complexes $\Omega_n$ it
is convenient to extend slightly the definitions and introduce a
class of more general cycle-free chessboard complexes. In
Section~\ref{sec:general} we introduce an even larger class of
hybrid chessboard complexes which contain $\Omega_n$ and
$\Delta_n$, as well as $\Delta(A)$ for $A\subset \mathbb{Z}^2$, as
special cases.

If $X$ is a finite set than $\Omega_X\subset \Delta_X=\Delta_{X,
X}$ is defined as the union of all simplices
$S_\rho=\{(x_1,x_2),(x_2,x_3),\ldots,(x_{n-1},x_n)\}$ where $\rho
: x_1\prec x_2\prec \ldots\prec x_n$ is a linear order on $X$.
More generally, given a bijection $\alpha : Y\rightarrow X$ of two
finite sets, let $\Omega(X\times Y;\alpha)$ be the complex of all
non-taking rook configurations in $X\times Y$ without
sub-configurations of the form $\{(x_1,y_2), (x_2,y_3),\ldots,
(x_m,y_1)\}$ where $1\leq m\leq\vert X\vert=n$ and
$x_{j}=\alpha(y_j)$ for each $j$. It is clear that all these
complexes are isomorphic to $\Omega_n$ if $\vert X\vert =\vert
Y\vert =n$.

\medskip
For visualization and as a convenient bookkeeping device,
simplices in $\Delta(X\times Y)$ as well as in $\Omega(X\times
Y;\alpha)$ can be represented as matchings in the complete
bipartite graph $K_{X,Y}$.

\begin{figure}[hbt]
\centering
\includegraphics[scale=0.60]{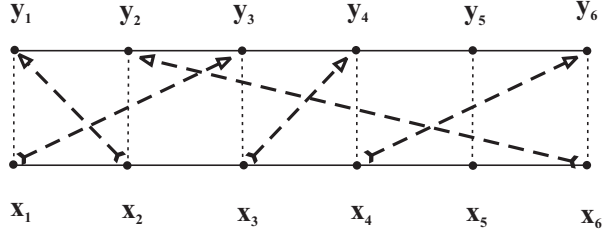}
\caption{A  cycle in $\Delta_6$.}\label{fig:bipartite}
\end{figure}
The partial matching
$\{(x_1,y_3),(x_2,y_1),(x_3,y_4),(x_4,y_6),(x_6,y_2)\}$, exhibited
in Figure~\ref{fig:bipartite}, clearly determines a non-taking
rook placement on the chessboard $X\times Y$ where
$X=\{x_i\}_{i=1}^6$ and $Y=\{y_i\}_{i=1}^6$. If $\alpha :
Y\rightarrow X$ is the bijection $y_j\mapsto x_j$ then this
matching does not contribute a simplex to $\Omega(X\times
Y;\alpha)$ since it contains a cycle
$$x_1\mapsto y_3\downarrow x_3\mapsto y_4\downarrow x_4\mapsto y_6
\downarrow x_6\mapsto y_2\downarrow x_2\mapsto y_1.$$

The following proposition  establishes a key structural property
for cycle-free chessboard complexes $\Omega_n$.

\begin{prop}\label{prop:link}
The link ${\rm Link}(v)={\rm Link}_{\Omega_n}(v)$ of each vertex
$v$ in the cycle-free chessboard complex $\Omega_n$ is isomorphic
to $\Omega_{n-1}$.
\end{prop}

\medskip\noindent
{\bf Proof:} Let us choose $\Omega(X\times Y;\alpha)$ as our model
for $\Omega_n$ where $X=\{x_i\}_{i=1}^n$ and $Y=\{y_i\}_{i=1}^n$
while the bijection $\alpha : Y\rightarrow X$ maps $y_j$ to $x_j$.
Let  $v=(x_i,y_j)\in X\times Y$, where $i\neq j$.

Define $X':=X\setminus\{x_i\}$ and $Y':=Y\setminus \{y_j\}$. Let
$\alpha' : Y'\rightarrow X'$ be the bijection defined by
$\alpha'(y_i)= x_j$ and $\alpha(y_k)=x_k$ for $k\in
[n]\setminus\{i,j\}$. Than it is not difficult to show that ${\rm
Link}_{\Omega_n}(v)\cong \Omega(X'\times Y',\alpha')\cong
\Omega_{n-1}$. \hfill$\square$

\subsection{$\Omega_n$ as a digraph complex}
\label{sec:digraph}

Chessboard complexes $\Delta_n$ and cycle-free chessboard
complexes $\Omega_n$, as well as their natural generalizations,
admit another, equally useful description as directed graph
(digraph) complexes.

A chessboard $A_n=[n]\times [n]$ is naturally interpreted as a
complete digraph $DK_n$ (loops included) where each $(i,j)\in A_n$
contributes a directed edge $\overrightarrow{ij}$ in $DK_n$. A
directed subgraph $\Gamma\subset DK_n$ describes an admissible
rook configuration in $A_n$ if and only if
\begin{figure}[hbt]
\centering
\includegraphics[scale=0.40]{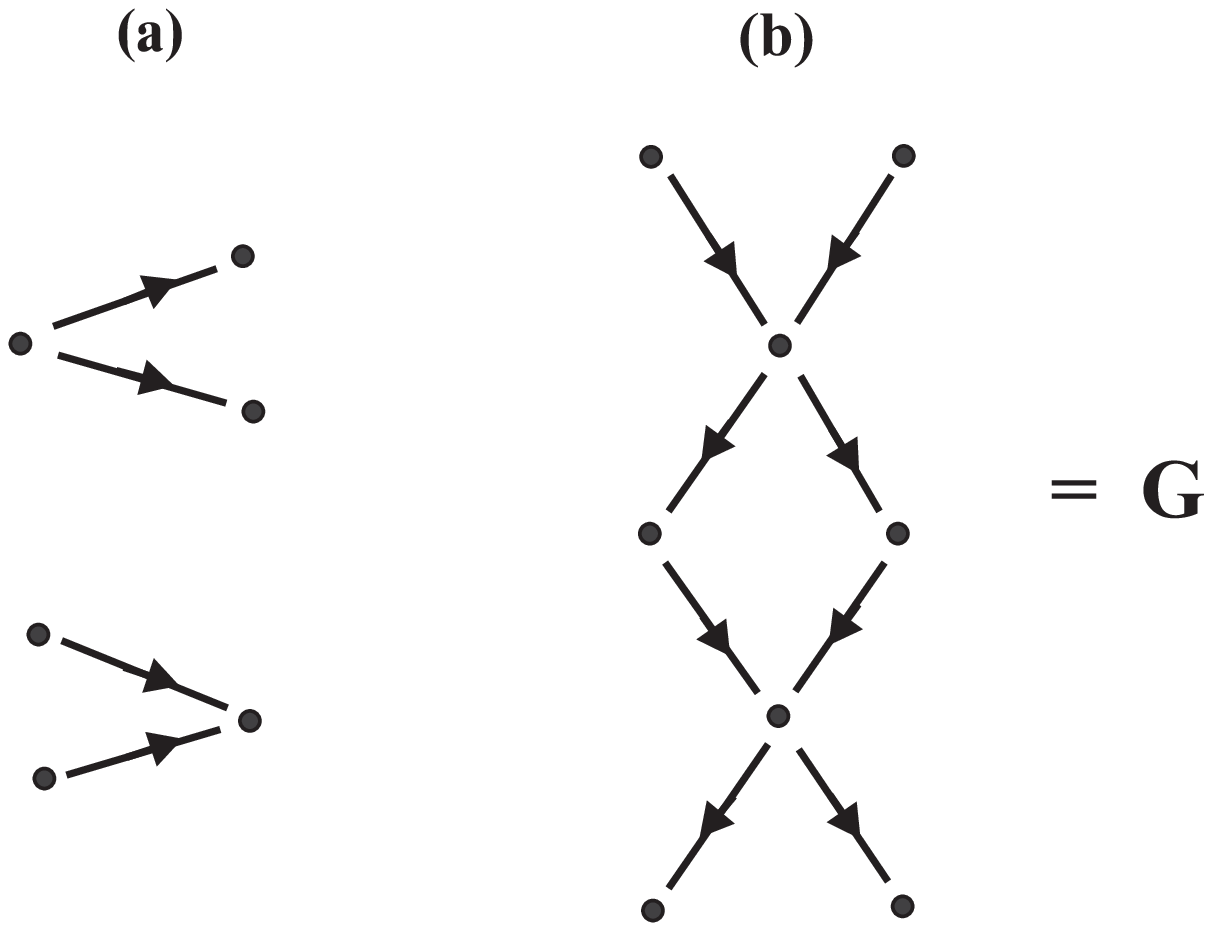}
\caption{\mbox{$\Omega(G)= S^0\ast S^0\ast S^0\ast S^0 = S^3$.}}
\label{fig:omega-1}
\end{figure}
no two directed edges in $\Gamma$ are allowed to have the same
{\em tail} or the same {\em end}. In other words configurations
depicted in Figure~\ref{fig:omega-1}~(a) are banned from the graph
$\Gamma$. It follows that $\Delta_n$ is the complex of all
subgraphs of $DK_n$ such that the associated connected components
are either directed cycles or directed paths. The complex
$\Omega_n$ arises as the cycle-free subcomplex of $\Delta_n$,
i.e.\ $\Gamma\in \Omega_n$ if only directed paths are allowed as
connected components of $\Gamma$. This definition reveals that
probably the closest relative of $\Omega_n$, that has been
systematically analyzed so far, is the complex $\Delta_n^{DM}$ of
directed matchings on the node set [n] introduced in
\cite{Bj-We-99}.

More generally, for each directed graph $G$ one can define the
associated complexes $\Delta(G)$ and $\Omega(G)$ as the complexes
of all directed subgraphs $\Gamma$ in $G$ which have only directed
paths and cycles (respectively paths alone) as connected
components. For example if $G$ is the directed graph depicted in
Figure~\ref{fig:omega-1}~(b) then $\Delta(G)=\Omega(G)\cong S^3$.

\section{Generalized cycle-free complexes}
\label{sec:general}

Let $\Omega(X\times Y,\alpha)$ be the cycle-free chessboard
complex associated to sets $X,Y\subset \mathbb{Z}$ and a bijection
$\alpha : Y\rightarrow X$. Assume that $A\subset \mathbb{Z}^2$ is
a finite superset of $X\times Y$. Define $\Omega =
\Omega(A,X\times Y,\alpha)$ as the subcomplex of the full
chessboard complex $\Delta(A)$ by the condition that $S\in
\Delta(A)$ is in $\Omega$ if and only if the restriction of $S$ on
$\Delta(X\times Y)$ is in $\Omega(X\times Y,\alpha)$. $\Omega$ is
referred to as the generalized cycle-free chessboard complex.

If $A = (X\cup Z)\times (Y\cup T)$, where $X\cap Z=\emptyset=Y\cap
T$, let $$\Omega^{Y,T}_{X,Z}:=\Omega(A,X\times Y,\alpha).$$ The
isomorphism type of the complex $\Omega^{Y,T}_{X,Z}$ depends only
on cardinalities of sets $X,Y,Z,T$ so if $\vert X\vert=\vert
Y\vert =n, \vert Z\vert =m,$ and $\vert T\vert=p$, we will
frequently denote by $\Omega^{n,p}_{n,m}$ one of its unspecified
representatives. If $p=0$ we write $\Omega_{n,m}:=
\Omega_{n,m}^{n,0}$ and if  $m=0$, the complex
$\Omega_{n,0}=\Omega_n$ reduces to the standard cycle-free
chessboard complex defined on a $n\times n$-chessboard.

\begin{defin}\label{def:reduced}
Let $\Omega = \Omega(A,X\times Y,\alpha)$ be a generalized,
cycle-free chessboard complex based on a chessboard $A\subset
\mathbb{Z}^2$, where $X\times Y\subset A$ and $\alpha :
Y\rightarrow X$ is an associated bijection. Let $v = (a,b)\in A$.
The $v$-{\em reduced complex} $\Omega' = \Omega'_v =
\Omega(A',X'\times Y',\alpha')$ of $\Omega$ is defined as follows.
Let $A':=A\setminus (\{a\}\times\mathbb{Z}\cup
\mathbb{Z}\times\{b\})$.

\begin{enumerate}
 \item[{\rm (a)}] If both $a\in X$ and $b\in Y$ let $X':=
X\setminus\{a\}, Y':=Y\setminus\{b\}$ and let $\alpha' :
Y'\rightarrow X'$ be the bijection defined by
$\alpha'(\alpha^{-1}(a)):=\alpha(b)$, and $\alpha'(z)=\alpha(z)$
for  $z\neq \alpha^{-1}(a)$.
 \item[{\rm (b)}] If $a\in X$ and $b\notin Y$ let $X':=
 X\setminus\{a\}, Y':=Y\setminus\{\alpha^{-1}(a)\}$ and $\alpha' : Y'\rightarrow
 X'$ is the restriction of $\alpha$ on $Y'$.
 \item[{\rm (c)}] If $b\in Y$ and $a\notin X$ let
$Y':=Y\setminus\{b\}, X':=X\setminus\{\alpha(b)\}$ and $\alpha' :
Y'\rightarrow  X'$ is the restriction of $\alpha$ on $Y'$.
 \item[{\rm (d)}] If neither $a\in X$ nor $b\in Y$, let $X'=X,
 Y'=Y$ and $\alpha'=\alpha$.
\end{enumerate}

\end{defin}

\medskip
The following proposition records for the future reference the key
structural property of generalized cycle-free chessboard complexes
$\Omega = \Omega(A,X\times Y,\alpha)$. The proof is similar to the
proof of Proposition~\ref{prop:link} so we omit the details.

\begin{prop}\label{prop:structural}
If ${\rm Link}(v)={\rm Link}_\Omega(v)$ is the link of a vertex
$v=(a,b)\in A$ in $\Omega=\Omega(A,X\times Y,\alpha)$ then there
is an isomorphism
$$
{\rm Link(v)\cong \Omega(A',X'\times Y',\alpha')}
$$
where $\Omega(A',X'\times Y',\alpha')$ is the $v$-reduced complex
of the generalized cycle-free chessboard complex $
\Omega(A,X\times Y,\alpha)$ (Definition~\ref{def:reduced}).
\end{prop}

\section{Filtrations of chessboard complexes}
\label{sec:filtrations}

The chessboard complex $\Delta(A)$ functorially depends on the
chessboard $A\subset \mathbb{Z}^2$. It follows that a filtration
$$
A_0\subset A_1\subset\ldots \subset A_{m-1}\subset A_m\subset A
$$
induces a filtration of the complex $\Delta(A)$,
$$
\Delta(A_0)\subset \Delta(A_1)\subset\ldots \subset
\Delta(A_{m-1})\subset \Delta(A_m)\subset \Delta(A).
$$
This filtration in turn induces a filtration
$\{F_j(\Omega)\}_{j=0}^m$ of the associated generalized,
cycle-free chessboard complex $\Omega=\Omega(A,X\times Y,\alpha)$.
If $X\times Y\subset A_0$ then clearly
$F_j(\Omega)=\Omega(A_j,X\times Y,\alpha)$. We are particularly
interested in filtrations where $A_j\setminus A_{j-1}=\{a_j\}$ is
a singleton. Consequently a filtration is determined once we
choose a linear order of the elements (elementary squares) of the
set $A\setminus A_0$.

\medskip
 A basic fact and a well known consequence of the {\em
Gluing Lemma} \cite{Brown} is that the homotopy type of the
``double mapping cylinder'' (homotopy colimit) of the diagram
$B\stackrel{f}\longleftarrow A \stackrel{g}\longrightarrow C$ of
spaces (complexes) depends only on homotopy types of maps $f$ and
$g$. It follows that if both maps $f$ and $g$ are homotopic to
constant maps the associated double mapping cylinder has the
homotopy type of a wedge $B\vee \Sigma(A)\vee C$. From here we
immediately deduce that if a simplicial complex $X = X_1\cup X_2$
is expressed as a union of its sub-complexes such that both $X_1$
and $X_2$ have the homotopy type of a wedge of $n$-dimensional
spheres while the intersection $X_1\cap X_2$ is a wedge of
$(n-1)$-dimensional spheres, then the complex $X$ is also a wedge
of $n$-dimensional spheres. An immediate consequence is the
following lemma.

\begin{lema}\label{lem:vertex}
Let $K$ be a finite simplicial complex. Given a vertex $v\in K$,
let ${\rm Link}_K(v)$ and ${\rm Star}_K(v)$ be the link and star
subcomplex of $K$. Let $A$-${\rm Star}_K(v) = K\setminus \{v\}$ be
the ``anti-star'' of $v$ in $K$, i.e.\ the complex obtained by
deleting $v$ from all simplices, or equivalently by removing the
``open star'' of $v$ from $K$. If $A$-${\rm Star}(v)$ is homotopy
equivalent to a wedge of $n$-dimensional spheres and ${\rm
Link}_K(v)$ is homotopy equivalent to a wedge of
$(n-1)$-dimensional spheres, then the complex $K$ itself has the
homotopy type of a wedge of $n$-dimensional spheres.
\end{lema}

One way of proving that a simplicial complex is homotopically a
wedge of $n$-spheres is to iterate Lemma~\ref{lem:vertex}. In the
following section we show that among the complexes where this
strategy can be successfully carried on are some generalized
cycle-free complexes.

\section{Complexes $\Omega_{n,m}$}

\begin{prop}\label{prop:primena}
The complex $\Omega_{n,m}$ is homotopy equivalent to a wedge of
$(n-1)$-dimensional spheres provided $m\geq n$.
\end{prop}

\medskip\noindent
{\bf Proof:} Let us establish the statement for all complexes
$\Omega_{n,m}$, where $m\geq n$, by induction on $n$. Note that
$\Omega_{2,2}$ is a circle and that $\Omega_{2,m}$ for $m\geq 3$
is always a connected, $1$-dimensional complex, hence a wedge of
$1$-spheres.

Assume, as an inductive hypothesis, that $\Omega_{n,m}$ is
homotopic to a wedge of $(n-1)$-spheres for each $m\geq n$.

Our model for $\Omega_{n,n}$ will be the complex $\Omega_{X,Z}^Y$
where $X=\{x_i\}_{i=1}^n, Y=\{y_i\}_{i=1}^n, Z=\{z_i\}_{i=1}^{n}$,
where $\alpha : Y\rightarrow X$ is the canonical bijection
$y_j\mapsto x_j$.

Our model for $\Omega_{n+1,n+1}$ will be the complex
$\Omega_{X',Z'}^{Y'}$ where $X'=X\cup\{x_0\}, Y'=Y\cup\{y_0\},
Z'=Z\cup\{z_0\}$ and the bijection $\alpha' : Y'\rightarrow X'$ is
the (unique) extension of $\alpha$ characterized by
$\alpha'(y_0)=x_0$.

\begin{figure}[hbt]
\centering
\includegraphics[scale=0.70]{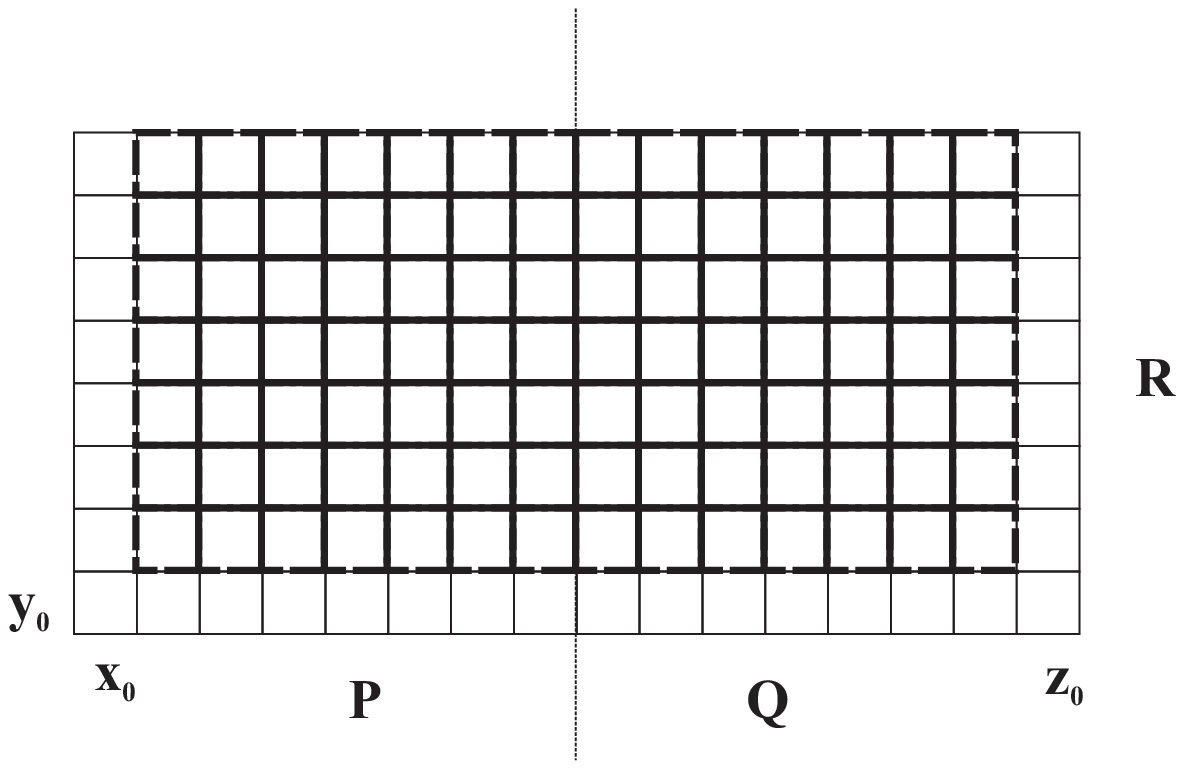}
\caption{}\label{fig:chess-1}
\end{figure}

Following the strategy outlined in Section~\ref{sec:filtrations},
we define a filtration of the complex $\Omega_{n+1,n+1}$ by
choosing $A_0=\{(z_0,y_0)\}\cup ((X'\cup Z)\times Y)$ as the
initial chessboard and selecting a linear order on the set
$W:=((X\cup Z)\times \{y_0\}) \cup (\{z_0\}\times Y)$ of
elementary squares (Figure~\ref{fig:chess-1}). Note that the
element $(x_0,y_0)$ is omitted since it is not allowed to be a
vertex of the cycle-free complex $\Omega(X',Y')$. Let
$$
 P = \{(x_i,y_0)\}_{i=1}^n, \quad
Q=\{(z_i,y_0)\}_{i=1}^{n}, \quad R=\{((z_0,y_i))\}_{i=1}^n.
$$
List elements of $W= P\cup Q\cup R$ in the order of appearance in
this $\cup$-decomposition. Within each of the blocks $P, Q, R$ the
elements can be ordered in an arbitrary way, say according to the
index $i=1,\ldots,n$.

If $W=\{v_k\}_{k=1}^{N}$ where $N=3n$, let
\begin{equation}
\{(z_0,y_0)\}\cup ((X'\cup Z)\times Y)=A_0\subset A_1\subset\ldots
\subset A_{N}=A=(X'\cup Z')\times Y' \label{eqn:filtration}
\end{equation}
be the filtration defined by $A_j:=A_0\cup\{v_k\}_{k=1}^j$. Let
$\{\Delta(A_j)\}_{j=0}^N$ be the associated filtration of the
chessboard complex $\Delta(A)$ and let $\{F_j(\Omega)\}_{j=0}^N$
be the induced filtration on the generalized cycle-free complex
$\Omega = \Omega(A,X'\times Y',\alpha')$. Note that
$F_j(\Omega)=\Omega(A_j,X'\times Y',\alpha)$ for $j\geq n$ while
in general $F_j(\Omega)=\Omega(A,X'\times Y',\alpha')\cap
\Delta(A_j)$.

By Proposition~\ref{prop:structural}, the homotopy type of the
link ${\rm Link}_k(v_k)$ of $v_k$ in the complex $F_k(\Omega)$ can
be described as follows.
\begin{enumerate}
\item[(I)]  $v_k\in P$, i.e.\ $v_k = (x_i,y_0)$ for some
$i=1,\ldots, n$.
$$
{\rm Link}_{k}(v_k) \cong \Omega_{n,n}
$$
\item[(II)]  $v_k\in Q$, i.e.\ $v_k = (z_i,y_0)$ for some
$i=1,\ldots, n$.
$$
{\rm Link}_{k}(v_k) \cong \Omega_{n,n}.
$$
\item[(III)] $v_k\in R$, i.e.\ $v_k=(z_0,y_i)$ for some
$i=1,\ldots, n$.
$$
{\rm Link}_{k}(v_k) \cong \Omega_{n,n+1}.
$$
\end{enumerate}
The complex $F_0(\Omega)$ is a cone with apex $(z_0,y_0)$, hence
it is contractible. In all cases (I)--(III), by the inductive
hypothesis, the complexes $\Omega_{n,n}$ and $\Omega_{n,n+1}$ have
the homotopy type of a wedge of $(n-1)$-dimensional spheres.
Consequently, by repeated use of Lemma~\ref{lem:vertex},
$\Omega_{n+1,n+1}$ has the homotopy type of a wedge of
$n$-dimensional spheres.

\medskip
It remains to be shown that the complex $\Omega_{n+1,m}$ has the
homotopy type of a wedge of $n$-dimensional spheres if $m>n+1$.
This is achieved by expanding the filtration
(\ref{eqn:filtration}) by adding vertices from new columns, in
some order, and applying the same argument as above.
\hfill$\square$

\section{Complexes $\Omega_{n,m}$ and the nerve lemma}

A classical result of topological combinatorics is the Nerve
Lemma. It was originally proved by J.~Leray in \cite{Le45}, see
also \cite{Bjo} for a more recent overview of applications and
related results.
\begin{lema}
{\bf (Nerve Lemma, \cite{Le45})} Let $\Delta$ be a simplicial
complex and $\{ L_i\}_{i=1}^k$ a family of subcomplexes such that
$\Delta = \cup_{i=1}^k~L_i.$ Suppose that every nonempty
intersection $L_{i_1} \cap L_{i_2} \cap \ldots \cap L_{i_t}$  is
$(\mu-t+1)$-connected for $t \geq 1.$ Then $\Delta$ is
$\mu$-connected if and only if ${\cal N}(\{L_i \}_{i=1}^k),$ the
nerve of the covering $\{L_i \}_{i=1}^k,$ is $\mu$-connected.
\end{lema}
 In the preceding section we showed that for $m\geq n$ the
complex $\Omega_{n,m}$ is a wedge of $(n-1)$-dimensional spheres,
consequently it is $(n-2)$-connected. Here we continue the
analysis of these complexes and establish a lower bound for the
connectivity of the complex $\Omega_{n,m}$ for any $m\geq 1$.
\begin{prop}
\label{prop:pomoc} The complex $\Omega_{n,m}$ is
$\mu_{n,m}$-connected, where $$\mu_{n,m}=\min \left\{\left[ \frac
{2n+m}3 \right]-2,n-2\right\}.$$
\end{prop}

\medskip\noindent
{\bf Proof:} We proceed by induction on $n$. For $n=2$, the
complex $\Omega_{2,1}$ is the union of two segments and so
non-empty (or $(-1)$-connected), and for $m\geq 2$ the complex
$\Omega_{2,m}$ is clearly connected (or $0$-connected).

Let us suppose that complexes $\Omega_{r,m}$ are
$\mu_{r,m}$-connected, whenever $r\leq n-1$, and consider the
complex $\Omega_{n,m}$. If $m\geq n$, then $\mu_{n,m}=n-2$, and
the complex $\Omega_{n,m}$ is $(n-2)$-connected by  Proposition
\ref{prop:primena}. Suppose that $1\leq m\leq n-1$, which implies
that $\mu_{n,m}\leq n-3$.

We use $\Omega_{[n],Z}^{[n],\emptyset}$, where $\vert Z\vert =m$,
as a model for the complex $\Omega_{n,m}$. For example, in order
to keep our chessboards in $\mathbb{Z}^2$, we could take
$Z=\{-1,-2,...,-m\}$. Let $\mathcal{L}_{n,m}=\{L_{z,i} \mid z\in
Z, i\in [n]\}$ be the family of subcomplexes of $\Omega_{n,m}$
where by definition $L_{z,i}:={\rm Star((z,i))}$ is the union of
all simplices with $(z,i)$ as a vertex, together with their faces.
Every maximal simplex in $\Omega_{n,m}$ must have a vertex
belonging to $Z\times [n]$. So, the collection $\mathcal{L}_{n,m}$
of contractible complexes is a covering of $\Omega_{n,m}$.

Let us apply the Nerve Lemma. It is easy to see that the
intersections of any $n-1$ complexes $L_{z,i}$ is nonempty. It
follows that the nerve $\mathcal{N}(\mathcal{L}_{n,m})$ of the
covering contains the full $(n-2)$-dimensional skeleton, hence it
is at least $(n-3)$-connected. It remains to show that the
intersection of any subcollection of $t$ of these complexes is at
least $(\mu_{n,m}-t+1)$-connected.

For the reader's convenience, we begin with the simplest case
$t=2$. There are three possibilities for the intersection
$L_{z_1,i}\cap L_{z_2,j}$.

\begin{itemize}

\item If $z_1\neq z_2$ and $i\neq j$, this intersection is a join
of the interval spanned by  vertices $(z_1,i),(z_2,j)$, and a
subcomplex of type $\Omega_{n-2,m}$. Therefore, it is
contractible.

\item If $z_1\neq z_2$ and $i=j$, this intersection is the
subcomplex of type $\Omega_{n-1,m-1}$, which is at least
$\mu_{n-1,m-1}=(\mu_{n,m}-1)$-connected by the induction
hypothesis.

\item If $z_1=z_2$ and $i\neq j$, this intersection is the
subcomplex of the type $\Omega_{n-2,m+1}$ which is
$\mu_{n-2,m+1}$-connected.

Then $\left[ \frac{2(n-2)+(m+1)}3\right]-2=\mu_{n,m}-1$. Also,
$(n-2)-2\geq \mu_{n,m}-1$ because $\mu_{n,m}\leq n-3$. Therefore,
$\mu_{n-2,m+1}\geq \mu_{n,m}-1$.
\end{itemize}

Similar arguments apply also in the case $t\geq 3$. The
intersection $L_{z_1,i_1}\cap L_{z_2,i_2}\cap \cdots \cap
L_{z_t,i_t}$ could be either contractible (when for some $h\in
\{1,2,...,t\}$ both $z_h$ and $i_h$ are different from all other
$z_j$ and $i_j$ respectively), or it could be a subcomplex of the
type $\Omega_{r,s}$ where both $r\geq n-t$ and $r+s\geq n+m-t$.
Then $2r+s\geq 2n+m-2t\geq 2n+m-3t+3$. Actually it could be easily
proved more, i.e. that $2r+s\geq 2n+m-\frac 32 t$, but we need the
more precise estimate only in the case $t=2$.

The above inequality implies $\left[ \frac {2r+s}3\right] -2\geq
\mu_{n,m}-t+1$. Also, $r-2\geq n-t-2\geq \mu_{n,m}-t+1$, because
$\mu_{n,m}\leq n-3$.

These two facts together imply that $\mu_{r,s}=\min \left\{ \left[
\frac {2r+s}3\right]-2,r-2\right\}\geq \mu_{n,m}-t+1$ which is
precisely the desired inequality. \hfill $\square$

\section{Complexes $\Omega_n$}

Now we are ready to prove our main result, i.e. to establish
high-connectivity of the complex $\Omega_n$.

\begin{prop}
For each $n\geq 5$, $\pi_1(\Omega_n)=0$.
\end{prop}

\medskip\noindent
{\bf Proof:} We apply the Nerve Lemma on the complex $L := L_1\cup
L_2\cup L_3$ where $L_j$ is the subcomplex of $\Omega_n$ based on
the chessboard $[n]\times ([n]\setminus\{i\})$. In other words a
simplex $\sigma\in\Omega_n$ is in $L_i$ if and only if it doesn't
have a vertex of the type $(\cdot, i)$.

It is clear that the $1$-skeleton of $\Omega_n$ is a subcomplex of
$L$, hence it suffices to show that $L$ is $1$-connected. Since
$L_1\cap L_2\cap L_2\neq\emptyset$ it is sufficient to show that
$L_i$ is $1$-connected for each $i$ and that $L_i\cap L_j$ is
connected for each pair $i\neq j$. Since $L_i\cong \Omega_{n-1,1}$
and $n\geq 5$ the first part follows from
Proposition~\ref{prop:pomoc}. Similarly, since $L_i\cap L_j\cong
\Omega_{n-2,2}$, again by Proposition~\ref{prop:pomoc} the complex
$L_i\cap L_j$ is connected if $n\geq 5$. \hfill $\square$

\begin{theo}
\label{thm:main} The complex $\Omega_n$ is $\mu_n$-connected,
where $\mu_n=\left[ \frac {2n-1}3\right] -2$.
\end{theo}

\medskip\noindent
{\bf Proof:}  For $n=2$ the complex $\Omega_2$ consists of two
points and is nonempty or $(-1)$-connected. For $n=3$ the complex
$\Omega_3$ is also nonempty ($(-1)$-connected), being an union of
two disjoint circles. The complex $\Omega_4$ is $0$-connected.
Indeed, each pair $v_0, v_1$ of vertices in $\Omega_4$ belongs to
a subcomplex isomorphic to $\Omega_{2,2}$ which is connected.

Let us assume that $n\geq 5$. We already know that
$\pi_1(\Omega_n)=0$ so it remains to be shown that
$H_j(\Omega_n)\cong 0$ for $j\leq \mu_n$. We establish this fact
by induction on $n$.

Let us suppose that the statement of the theorem is true for
complexes $\Omega_{n-2}$ and $\Omega_{n-1}$. Consider the
subcomplex $\Theta_n$ of $\Omega_n$ formed by simplices having
possibly a vertex of the type $(1,i)$ or $(j,1)$ but not both.
Here is an excerpt from the long homology exact sequence of the
pair $(\Omega_n,\Theta_n)$.
\begin{equation}
\label{long} \cdots \to H_{\mu_n}(\Theta_n)\to
H_{\mu_n}(\Omega_n)\to H_{\mu_n}(\Omega_n,\Theta_n)\to \cdots
\end{equation}

We need yet another exact sequence involving  complexes $\Omega_n$
and $\Theta_n$. For motivation, the reader is referred to
\cite{ShaWa04} where  similar sequences are constructed in the
context of usual chessboard complexes.

Let us denote by $\Theta_n^1$ the subcomplex of $\Theta_n$
consisting of simplices having one vertex of the type $(1,i)$, and
by $\Theta_n^2$ the subcomplex of $\Theta_n$ consisting of
simplices having one vertex of the type $(j,1)$. We use the
Mayer-Vietoris sequence for the decomposition
$\Theta_n=\Theta_n^1\cup \Theta_n^2$. Obviously $\Theta_n^1\cap
\Theta_n^2=\Omega_{n-1}$ so we obtain the following exact sequence
\begin{equation}
\label{Mayer} \cdots \to H_{\mu_n}(\Theta_n^1)\oplus
H_{\mu_n}(\Theta_n^2) \to H_{\mu_n}(\Theta_n) \to
H_{\mu_n-1}(\Omega_{n-1}) \to \cdots
\end{equation}

\noindent Since both $\Theta_n^1$ and $\Theta_n^2$ are the
complexes of type $\Omega_{n-1,1}$, they are
$\mu_{n-1,1}$-connected by Proposition \ref{prop:pomoc}. Since
$\mu_{n-1,1}=\left[ \frac {2n-1}3\right]-2=\mu_n$ we observe that
$H_{\mu_n}(\Theta_n^1)\oplus H_{\mu_n}(\Theta_n^2)=0$.

The complex $\Omega_{n-1}$ is $\mu_{n-1}$-connected by the
induction hypothesis, and $\mu_{n-1}=\left[ \frac {2n-3}3\right]
-2\geq \left[ \frac {2n-1}3\right] -2-1=\mu_n-1$. Therefore,
$H_{\mu_n-1}(\Omega_{n-1})=0$.

These facts, together with the exactness of the sequence
(\ref{Mayer}), allow us to conclude that $H_{\mu_n}(\Theta_n)=0$.

The homology of the pair $(\Omega_n,\Theta_n)$ is isomorphic to
the homology of the quotient $\Omega_n/\Theta_n$. If we denote by
$I_{i,j}$ (for $i\neq j$) the $1$-simplex with endpoints $(1,i)$
and $(j,1)$, the argument similar to the one from Proposition
\ref{prop:link} shows that this quotient is homotopy equivalent to
the wedge
$$\bigvee_{1<i\neq j\leq n} (I_{i,j}\ast \Omega_{n-2})/(\partial
I_{i,j}\ast \Omega_{n-2}).$$

Each quotient $(I_{i,j}\ast \Omega_{n-2})/(\partial I_{i,j}\ast
\Omega_{n-2})$  is homotopy equivalent to a wedge of double
suspensions of the complex $\Omega_{n-2}$. These double
suspensions are by the induction hypothesis
$(\mu_{n-2}+2)$-connected, and $(\mu_{n-2}+2)=\left[ \frac
{2n-5}3\right] \geq \left[ \frac {2n-1}3\right] -2=\mu_n$.
Therefore $H_{\mu_n}(\Omega_n,\Theta_n)=0$.

Finally, from the exact sequence (\ref{long}) we deduce
$H_{\mu_n}(\Omega_n)=0$, which completes our inductive argument.
\hfill $\square$
\bigskip

Substituting $n=p+1$ and taking the suspension, one immediately
obtains the desired estimate for the connectivity of the complex
$\mbox{Sym}^{(p)}_*$ introduced by Ault and Fiedorowicz in
\cite{AuFie07}.

\begin{cor}\label{cor:rezultat}
\label{posl} The complex $\mbox{Sym}^{(p)}_*$ is $\left[ \frac
23(p-1)\right]$-connected.
\end{cor}

\medskip\noindent
{\bf Proof:} Since by definition $\mbox{Sym}^{(p)}_*=\Sigma
\Omega_{p+1}$ it is $\gamma_p$-connected where
$$\gamma_p=\left[ \frac {2(p+1)-1}3\right] -2+1=\left[ \frac 23
(p-1)\right].$$ \hfill $\square$

\section{Tightness of the bound}
\label{sec:tightness}

Our objective in this section is to explore how far from being
tight is the connectivity bound established in
Theorem~\ref{thm:main}. Our central result is
Theorem~\ref{thm:tight1} which says that the constant $\mu_{n}$ is
the best possible at least if $n=3k+2$ for some $k\geq 1$.

\subsection{The case $n=3k+2$}

It is well known that $H_2(\Delta_{5,5})\cong \mathbb{Z}_3$,
\cite{BLVZ}, \cite{ShaWa04}, \cite{Jo07-2}. This fact was
essentially established in \cite{BLVZ}, Proposition~2.3.
Unfortunately the proof of this proposition suffers from an easily
rectifiable error which was detected too late to be inserted in
the final version of \cite{BLVZ}. Since the proof of
Proposition~\ref{prop:pet-puta-pet} depends on this result, we
start with a proposition which isolates the needed fact, points to
the error in the original proof of Proposition~2.3. and shows how
it should be corrected.

\begin{figure}[hbt]
\centering
\includegraphics[scale=0.5]{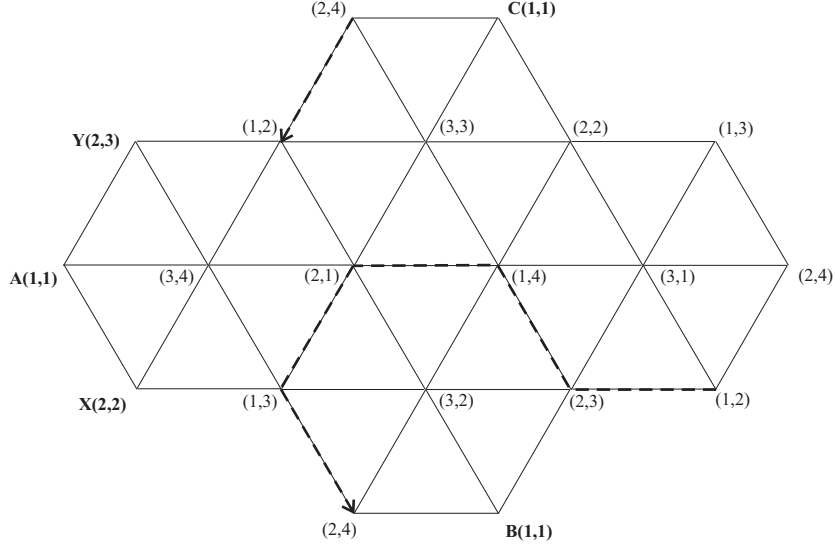}
\caption{The chessboard complex
$\Delta_{3,4}$.}\label{fig:torus-2}
\end{figure}

Recall that the chessboard complex $\Delta_{3,4}$ is isomorphic to
a torus $T^2$. More precisely, \cite{BLVZ}, p.\ 30, the universal
covering space of $\Delta_{3,4}$ is the triangulated honeycomb
tessellation of the plane. An associated fundamental domain for
$\Delta_{3,4}=\mathbb{R}/\Gamma$ is depicted in
Figure~\ref{fig:torus-2} with the lattice $\Gamma$ generated by
vectors $x=\overrightarrow{AB}$ and $y=\overrightarrow{AC}$. As
clear from the picture, $x=4a +2b$ and $y=2a+4b$ are generators of
the lattice $\Gamma:=H_1(\Delta_{3,4})$, where
$a:=\overrightarrow{AX}$ and $b:=\overrightarrow{AY}$. If
$\Gamma_1$ is the lattice spanned by vectors $6a$ and $6b$ then
$\Gamma_1\subset\Gamma$ and $\Gamma/\Gamma_1\cong \mathbb{Z}_3$.
As a consequence (\cite{BLVZ}, Lemma~2.2.),
$${\rm Coker}(H_1(\Delta_{3,3}) \rightarrow
H_1(\Delta_{3,4}))\cong \Gamma/\Gamma_1 \cong \mathbb{Z}_3.$$

\begin{prop}\label{prop:korekcija}
There is an isomorphism
$$
H_2(\Delta_{5,5}) \cong \oplus_{i=1}^4 H_1(\Delta^i_{3,4})/N\cong
\Gamma^{\oplus 4}/N \cong \mathbb{Z}_3
$$
where $\Delta^i_{3,4}\cong \Delta_{3,4}$ for each $i$ and $N=A+B$,
where $A=\Gamma_1^{\oplus 4}$ and $B = {\rm Ker}(\Gamma^{\oplus
4}\stackrel{\theta}\rightarrow \Gamma)$,
$\theta(x,y,z,t)=x+y+z+t$.
\end{prop}

\medskip\noindent
{\bf Proof:} The proof follows into the footsteps of the proof of
Proposition~2.3.\ from \cite{BLVZ}. The only defect in the proof
of that proposition is an incorrect determination of the kernel
${\rm Ker}(\gamma)$ of the homomorphism $\gamma :
H_1(\Delta_{3,3}^i)\rightarrow H_1(\overline{\Delta}_{4,3})$ in
the commutative diagram (loc.\ cit.), leading to the omission of
the factor $B$ in the decomposition $N=A+B$. As a consequence, the
group $H_2(\Delta_{5,5})$ is isomorphic to the group
$\Gamma^{\oplus 4}/N\cong \Gamma/\Gamma_1\cong \mathbb{Z}_3$,
rather than to group $\Gamma^{\oplus 4}/A\cong
\mathbb{Z}_3^{\oplus 4}$, as erroneously stated in the formulation
of Proposition~2.3.\ in \cite{BLVZ}. \hfill $\square$

\begin{figure}[hbt]
\centering
\includegraphics[scale=0.80]{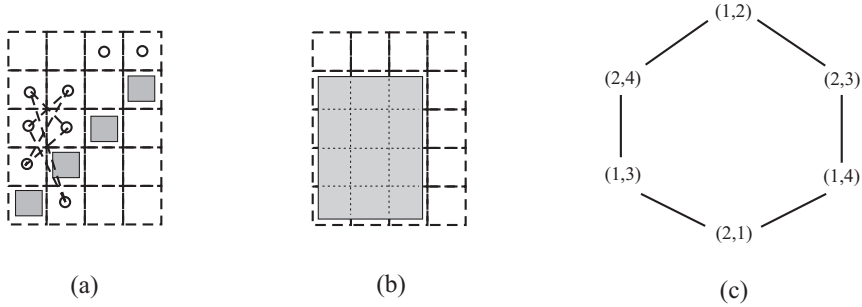}
\caption{A non-zero class in
$H_2(\Omega_5)$.}\label{fig:cetri-pet}
\end{figure}

\begin{exam}\label{exam:primena}{\rm
Proposition~\ref{prop:korekcija} is in practise applied as
follows. Suppose we want to check if a subcomplex $S\subset
\Delta_{4,5}\subset \Delta_5$ contributes a non-trivial
$2$-dimensional class to $H_2(\Delta_5)$. For example let $S\cong
S^1\ast S^0$ be the $2$-sphere shown in
Figure~\ref{fig:cetri-pet}~(a) where $S^1$ is the hexagon shown in
Figure~\ref{fig:cetri-pet}~(c) and $S^0=\{(3,5),(4,5)\}$. Let
$\Delta^i_{3,4},\, 1\leq i\leq 4$, be the chessboard complex
associated to the chessboard $[4]\setminus\{i\}\times [4]$ so for
example $\Delta_{3,4}^4\cong\Delta_{3,4}\cong\Gamma$ is associated
to the board depicted in Figure~\ref{fig:cetri-pet}~(b). Recall
that $\Delta_{4,5}/\Delta_{4,4}\cong
\vee_{i=1}^4~\Sigma(\Delta_{3,4}^i)$. Let $\nu :
H_2(\Delta_{4,5})\rightarrow \oplus_{i=1}^4
H_1(\Delta_{3,4}^i)\cong \Gamma^{\oplus 4}$ be a homomorphism
associated to the natural projection $\Delta_{4,5}\rightarrow
\Delta_{4,5}/\Delta_{4,4}$. Then the fundamental class $[S]$ is a
non-trivial element in $H_2(\Delta_5)$ if and only if $\nu([S])$
is not an element of $N=A+B$.

For example in our case the image of $[S]$ in $\Gamma^{\oplus 4}/N
\cong \Gamma/\Gamma_1\cong \mathbb{Z}_3$ is equal to the image of
the class of the circle $S^1$ depicted in
Figure~\ref{fig:cetri-pet}~(c) in $\Gamma/\Gamma_1$. By inspection
of Figure~\ref{fig:torus-2} we observe that this class is a
generator of $\Gamma$, hence $[S]$ is a generator of
$H_2(\Delta_{5})$. }
\end{exam}

\begin{prop}\label{prop:pet-puta-pet}
The inclusion $\Omega_5\hookrightarrow \Delta_{5}$ induces an
epimorphism $$H_2(\Omega_5)\stackrel{\alpha}{\longrightarrow}
H_2(\Delta_5)\cong \mathbb{Z}_3.$$ Moreover, for a class $[S]$
such that $\alpha([S])$ is a generator in $H_2(\Delta_5)$, one can
choose the fundamental class of the $2$-sphere $S\cong S^1\ast
S^0\cong \Omega_{2,2}\ast \Delta_{2,1}\subset \Omega_5$, depicted
in Figure~\ref{fig:cetri-pet}, where $\Omega_{2,2}\subset
\Delta_{[2],[4]}$ and $\Delta_{2,1}\cong \Delta(\{(3,5),(4,5)\})$.
\end{prop}

\medskip\noindent
{\bf Proof:} We have already demonstrated in
Example~\ref{exam:primena} that the image of $[S]$ in
$H_2(\Delta_5)$ is non-zero so the proof follows from the
observation that $S\subset \Omega_5$. \hfill $\square$

\begin{figure}[hbt]
\centering
\includegraphics[scale=0.80]{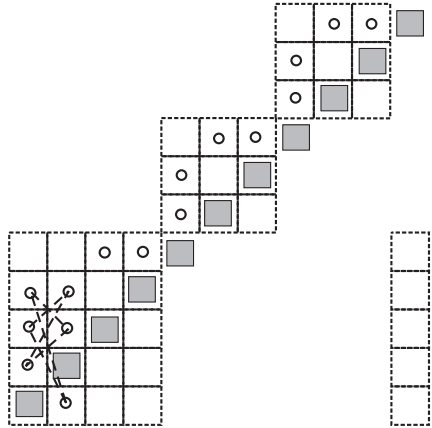}
\caption{The complex $\Omega_{11}$ is not
$6$-connected.}\label{fig:jedanaest}
\end{figure}

\begin{theo}\label{thm:tight1}
The inclusion map $\Omega_{3k+2}\hookrightarrow \Delta_{3k+2}$
induces a non-trivial homomorphism
$$
H_{2k}(\Omega_{3k+2})\longrightarrow H_{2k}(\Delta_{3k+2}).
$$
It follows that $H_{2k}(\Omega_{3k+2})$ is non-trivial, hence the
cycle-free chessboard complex $\Omega_{3k+2}$ is
$(2k-1)$-connected but not $(2k)$-connected for each $k\geq 1$.
\end{theo}

\medskip\noindent
{\bf Proof}: We already know that the result is true in the case
$k=1$. The general case is not much more difficult to prove in
light of the properties of chessboard complexes of the form
$\Delta_{3k+2}$ established in \cite{ShaWa04}. For example
Theorem~5.4.\ (loc.\ cit.) implies that
$H_{2k}(\Delta_{3k+2})\cong \mathbb{Z}_3$. Moreover, a generator
of this group is determined by a sphere $S^{2k}\cong
S^0\ast\ldots\ast S^0$ obtained as a join of $(2k+1)$ copies of
$S^0$ such that $(k+1)$ of them are vertical and the remaining $k$
are horizontal ``dominoes'', i.e.\ complexes of the form
$\Delta_{2,1}$ and $\Delta_{1,2}$ respectively. It is often
convenient to represent two dominoes of different type inside a
chessboard complex of the type $\Delta_{3,3}$, two of these
$(3\times 3)$-chessboards with pairs of complementary dominos are
indicated in Figure~\ref{fig:jedanaest}.

Let us illustrate the argument leading to the proof of the theorem
in the case of the complex $\Omega_{11}$, the proof of the general
case follows exactly the same pattern. Figure~\ref{fig:jedanaest}
exhibits a sphere $\Sigma:=S\ast S^1\ast S^1\cong S^6$, where $S$
is the $2$-sphere described in Example~\ref{exam:primena} while
the two copies of $S^1$ arise from the dominos in two $(3\times
3)$-blocks. It is clear that $\Sigma\subset \Omega_{11}$ so it
remains to be shown that the image of $\Sigma$ in $\Delta_{11}$
defines a non-zero homology class.

The image $\nu([S])$ of the class $[S]\in H_2(\Omega_5)$ in
$H_2(\Delta_5)$ is shown in Proposition~\ref{prop:pet-puta-pet} to
be non-trivial hence, according to Theorem~5.4.\ from
\cite{ShaWa04}, it is homologous (in $\Delta_5$) to a sphere
$S_1=S^0\ast S^0\ast S^0$ where two of the ``dominoes'' $S^0$ are
vertical. Hence $[\Sigma]$ is homologous (in $\Delta_{11}$) to the
fundamental class $[\Sigma_1]$ of $\Sigma_1:=S_1\ast S^1\ast S^1$
which, again by Theorem~5.4.\ from \cite{ShaWa04}, is non-trivial.
This completes the proof of the theorem. \hfill $\square$

\subsection{The cases $n=3k$ and $n=3k+1$}
\label{sec:cases}

Unfortunately the methods used in this paper do not alow us to
clarify if the constant $\mu_n$ from Theorem~\ref{thm:main} is the
best possible if $n=3k$ or $n=3k+1$ for some $k\geq 1$.
Nevertheless we are able to show that this bound should not be
expected to be too far off the actual bound.

\begin{prop}\label{prop:tight-1}
The group $H_{2k-1}(\Omega_{3k,1})$ is non-trivial. Moreover, a
non-trivial element of this group arises as the fundamental class
$\xi_{2k-1} = [\Sigma_{2k-1}]$ of a subcomplex
$\Sigma_{2k-1}\subset \Omega_{3k,1}$ isomorphic to the join
$S^0\ast\ldots\ast S^0\cong S^{2k-1}$ of\, $2k$ copies of\,
$0$-dimensional spheres.
\end{prop}

\begin{figure}[hbt]
\centering
\includegraphics[scale=0.80]{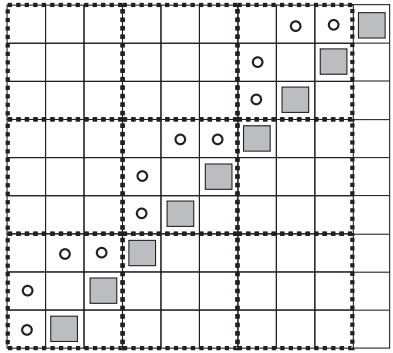}
\caption{A  cycle in $\Delta_9$.}\label{fig:tight-1}
\end{figure}

\medskip\noindent
{\bf Proof:} Our model for $\Omega_{n,1}=\Omega_{3k,1}$ is the
complex $\Omega_{X,Z}^{Y,\emptyset}$ where $X=\{x_1,\ldots,
x_n\}$, $Y=\{y_1,\ldots,y_n\}$, $Z=\{x_0\}$ and the bijection
$\alpha : Y\rightarrow X$ maps $y_j$ to $x_j$. The case $k=3$ is
depicted in Figure~\ref{fig:tight-1} where the shaded squares
correspond (from left to right) to squares $(x_j,y_j)$ while the
first column is filled with squares $(x_0,y_j),\, j=1,\ldots, n$.

Let $T_1:=\{(x_0,y_1), (x_0,y_2)\}, T_2:=\{(x_1,y_3),(x_2,y_3)\},
T_3:=\{(x_3,y_4),(x_3,y_5)\}, \ldots,$
$T_{2k-1}:=\{(x_{3k-3},y_{3k-2}),(x_{3k-3},y_{3k-1})\},
T_{2k}:=\{(x_{3k-2},y_{3k}),(x_{3k-1},y_{3k})\}$. Define
$\Sigma_{2k-1}$ as the join $T_1\ast\ldots\ast T_{2k}$. The proof
is completed by the observation that the cycle $\xi_{2k-1}$,
determined by the sphere $\Sigma_{2k-1}$, does not bound even in
the larger chessboard complex $\Delta_{X\cup\{x_0\},Y}$, cf.\
\cite{ShaWa04}, Section~3. \hfill $\square$

\medskip
The following corollary provides evidence that the connectivity
bound established in Theorem~\ref{thm:main} is either tight or
very close to the actual connectivity bound in the two remaining
cases, $n=3k, n=3k+1$.

\begin{cor} For each $k\geq 1$,
$$ \mbox{ {\rm either} }\quad H_{2k-1}(\Omega_{3k})\neq 0 \quad\mbox{ {\rm or} }
\quad H_{2k-1}(\Omega_{3k+1})\neq 0. $$
\end{cor}

\medskip\noindent
{\bf Proof:} Let $\Omega_{3k+1}$ be the cycle-free chessboard
complex based on the chessboard $[3k+1]\times [3k+1]$. Define
$\Omega_{3k,1}$ as the subcomplex of $\Omega_{3k+1}$ such that a
simplex $S\in \Omega_{3k+1}$ is in $\Omega_{3k,1}$ if and only if
$S\cap (\{1\}\times [3k+1])=\emptyset$. The quotient complex
$\Omega_{3k+1}/\Omega_{3k,1}$ has the homotopy type of a wedge
$\bigvee_{i=1}^{3k+1}\Sigma(\Omega_{3k}^{(i)})$ where each of the
complexes $\Omega_{3k}^{(i)}$ is isomorphic to $\Omega_{3k}$.
Consider the following fragment of the long exact sequence of the
pair $(\Omega_{3k+1},\Omega_{3k,1})$,
$$
\ldots\rightarrow\oplus_{i=1}^{3k+1} H_{2k-1}(\Omega_{3k}^{(i)})
\rightarrow H_{2k-1}(\Omega_{3k,1}) \rightarrow
H_{2k-1}(\Omega_{3k+1})\rightarrow\ldots
$$
The desired conclusion follows from the fact that
$H_{2k-1}(\Omega_{3k,1})\neq 0$. \hfill$\square$

\bigskip\noindent {\bf Conjecture}: The connectivity bound given in
Theorem~\ref{thm:main} is the best possible or in other words for
each $n\geq 2$,
$$H_{\mu_n+1}(\Omega_n)\neq 0 . $$

\section{Relatives of $\Omega_n$}
\label{sec:comparison}

The closest relative of $\Omega_n$, that has so far appeared in
the literature, is the complex $\Delta_n^{DM}$ of directed
matchings introduced by Bj\" orner and Welker in \cite{Bj-We-99},
see also Section~\ref{sec:digraph}. In this section we describe a
natural ``ecological niche'' for all these complexes and briefly
compare their connectivity properties.

In the sequel we put more emphasis on the directed graph
description of $\Omega_n, \Delta_n$ and related complexes
(Section~\ref{sec:digraph}). We silently identify a directed graph
with its set of directed edges (assuming the set of vertices is
fixed and clear from the context).

Let $DK_n$ be the complete directed graph on the set $[n]$ of
vertices (directed loops included) and $K_n^{{\uparrow}}$ its
companion with all loops excluded.

Following $\cite{Bj-We-99}$, let $\Delta_n^{DM}$ be the directed
graph complex of all {\em directed matchings} in $K_n^{\uparrow}$.
By definition, $\Gamma\subset K_n^{\uparrow}$ is a directed
matching if both the inn-degree and out-degree of each vertex is
at most one. This is equivalent to the condition that two graphs
depicted in Figure~\ref{fig:omega-1}~(a) are banned from $\Gamma$.
It follows that $\Gamma\subset DK_n$ is in $\Delta_n^{DM}$ if and
only if the connected components of $\Gamma$ are either directed
paths or directed cycles of length at least $2$, i.e.\ the only
difference between $\Delta_n$ and $\Delta_n^{DM}$ is that in the
former complex the cycles of length one (loops) are allowed.

Summarizing, if $\Gamma\subset DK_n$ then
\begin{enumerate}
 \item[(1)] $\Gamma\in\Delta_n \Leftrightarrow$ Each connected
component of $\Gamma$ is either a directed path or a directed
cycle,
 \item[(2)] $\Gamma\in\Delta_n^{DM} \Leftrightarrow$ Each connected
component of $\Gamma$ is either a directed path or a directed
cycle of length at least $2$,
 \item[(3)] $\Gamma\in\Omega_n \Leftrightarrow$ Each connected
 components of $\Gamma$ is a  directed path.
\end{enumerate}
The following definition introduces (some of) natural intermediate
complexes which interpolate between $\Omega_n$ and $\Delta_n$,
respectively $\Omega_n$ and $\Delta_n^{DM}$.

\begin{defin}\label{def:filtration}
Let $F_p^I=F_p(\Delta_n)$, respectively
$F_p^{II}=F_p(\Delta_n^{DM})$, be the subcomplex of $\Delta_n$
(respectively $\Delta_n^{DM}$) such that $\Gamma\in F_p(\Delta_n)$
(respectively $\Gamma\in F_p(\Delta_n^{DM})$) if and only if the
number of cycles, among the connected components in $\Gamma$ is at
most $p$.
\end{defin}

\begin{defin}
A graph $C\subset DK_n$ is a $p$-multicycle if $C=C_1\uplus
C_2\uplus\ldots\uplus C_p$ has exactly $p$ connected components
and each $C_j$ is a cycle. If $l_j:=l(C_j)$ is the length of $C_j$
then the multiset $t(C):=(l_1,l_2,\ldots,l_p)$ is called the type
of $C$ and the number $l(C):=l(C_1)+\ldots +l(C_k)$ is called the
length of the $p$-multicycle $C$. Let $\mathcal{C}_p$ be the set
of all $p$-multicycles and $\mathcal{C}_p^{\geqslant 2}$ the set
of all $p$-multicycles $C$ of type $t(C):=(l_1,l_2,\ldots,l_p)$
such that $l_j\geq 2$ for each $j$.
\end{defin}

\begin{prop}\label{prop:filtration}

 \begin{equation}\label{eqn:prva}
 F_0^I=F_0^{II}=F_0(\Delta_n)=F_0(\Delta_n^{DM})=\Omega_n
 \end{equation}
 \begin{equation}\label{eqn:druga}
 F_p^I/F_{p-1}^I\simeq \bigvee_{C\in \mathcal{C}_p} S^{l(C)-1}\ast
 \Omega_{n-l(C)}\cong \bigvee_{C\in \mathcal{C}_p} \Sigma^{l(C)}(\Omega_{n-l(C)})
 \end{equation}
 \begin{equation}\label{eqn:treca}
 F_p^{II}/F_{p-1}^{II}\simeq \bigvee_{C\in \mathcal{C}_p^{\geqslant 2}} S^{l(C)-1}\ast
 \Omega_{n-l(C)}\cong \bigvee_{C\in \mathcal{C}_p^{\geqslant 2}} \Sigma^{l(C)}(\Omega_{n-l(C)})
 \end{equation}
\end{prop}
The associated exact (spectral) sequences show that all these
complexes are closely related, in particular have very similar
connectivity properties. Let $\mu_n=[\frac{2n-1}{3}]-2$ and
$\nu_n=[\frac{2n+1}{3}]-2$. The complex $\Delta_n$ is
$\nu_n$-connected, as demonstrated by Bj\" orner et al.\ in
\cite{BLVZ}. The same connectivity bound was established by Bj\"
orner and Welker for $\Delta_n^{DM}$ in \cite{Bj-We-99}. Both
bounds are tight as proved by Shareshian and Wachs in
\cite{ShaWa04}. It follows from Proposition~\ref{prop:filtration}
that majority of complexes $F_p^I$ and $F_p^{II}$ share this
connectivity bound. On the other hand $\Omega_n$ is by
Theorem~\ref{thm:main} $\mu_n$-connected, hence all these
complexes are $\mu_n=\nu_n$ connected if $n=3k+2$ for some $k$.

\end{document}